\documentclass[10pt]{amsart}
\usepackage{amsfonts,amsmath,amssymb,amsthm}
\usepackage{t1enc}

\newtheorem{theorem}{Theorem}[section]
\newtheorem{lemma}[theorem]{Lemma}
\newtheorem{proposition}[theorem]{Proposition}
\newtheorem{corollary}[theorem]{Corollary}
\theoremstyle{definition}
\newtheorem{definition}[theorem]{Definition}

\newtheorem{remark}[theorem]{Remark}
\numberwithin{equation}{section}

\def\R{{\mathbb R}}
\def\im{{\text {\rm Im}}\,}
\def\re{{\text {\rm Re}}\,}
\def\ee{{\text {\rm e}}} 
\def\C{{\mathbb C}}

\def\dd{{\,{\text {\rm d}}}}

\def\ii{{\text {\rm i}}}

\def\ep{{\varepsilon}}
\def\de{{\delta}}
\def\ka{{\kappa}}
\newcommand{\comm}[2]{[ #1, #2]}

\newcommand{\abs}[1]{\vert{#1}\vert}

\newcommand{\opn}{\operatorname}

\renewcommand{\tilde}{\widetilde}

\begin{document}

\date{}

\title[Inverse scattering]{Inverse scattering for the 1-D Helmholtz equation}
\author{Ingrid Belti\c t\u a and Renata Bunoiu}
\address{Institute of Mathematics ``Simion Stoilow'' 
of the Romanian Academy, 
P.O. Box 1-764, Bucharest, Romania}
\address{Universit\' e de Lorraine, IECL, UMR 7502, 
Metz, France }
%\dedicatory{IBRB-vers.tex, 28.07.2014}

%\parskip=5pt

\begin{abstract}

We prove a uniqueness result for Nevanlinna functions. and  
this result is then used to give an elementary proof
of the  uniqueness in the inverse scattering problem for the equation
$ u'' + \frac{k^2}{c^2}u=0 $
on $\R$. Here $c$ is a real positive measurable function that is  bounded from below by a positive constant,
and is close to $1$ at $\pm \infty$. 
\end{abstract}

\maketitle

\section{Introduction}

The uniqueness for the inverse scattering problem for the equation 
\begin{equation}\label{bbw-1}
- u'' + q u = k^2 w u, 
\end{equation}
where $q\ge 0$ and $w$ is real, and  not required to have fixed sign,  was proved in \cite{BBW}.
Our motivation is to give a more elementary proof for a simplified problem, where 
from all the features that make \eqref{bbw-1} difficult, only the non-regularity of $w$ is retained.

In fact, 
we consider the inverse scattering problem for the equation
\begin{equation}\label{eq:4.1.1}
-u''-\frac{k^{2}}{c^{2}}u=0
\end{equation}
where $c$ is a real measurable function that satisfies
\begin{itemize}
\item[[H1]]  There exist $c_0>0$, $c_M>0$ such $c_{0}<c(x)\le c_{M}$, 
a.e.  $x\in \R$;
\item[[H2]]  The function $c-1$ belongs to $L^1(\R)$. 
%One has  $|c(x)-1|\le C(1+ \vert x\vert^2)^{-(1+\delta)/2}$, 
%a. e.  $x\in \R$, where  
%$\delta>0$, and $C$ is a positive constant independent of $x$. 
\end{itemize}

Then for each $k\in\R\setminus\{ 0 \}$ there exist unique solutions 
 $u_{1}(x,k)$, $u_{2}(x,k)$
to the Helmholtz equation \eqref{eq:4.1.1}
such that
$
u_{1}(x,k)\sim \ee^{\ii kx}$ when $x\to\infty$ and 
$u_{2}(x,k)\sim \ee^{-\ii kx}$ when $x\to -\infty$. 
Then $u_1(\cdot, k)$ and $\overline{u_1(\cdot,k)}$ 
(respectively $u_2(\cdot, k)$ and $\overline{u_2(\cdot,k)}$)
are linearly independent solutions of  \eqref{eq:4.1.1},
therefore
$$
\begin{aligned}
 u_{1}(x,k)\sim &\frac{1}{T_{2}(k)}
\ee^{\ii kx}+\frac{R_{2}(k)}{T_{2}(k)}\ee^{-\ii kx}
\quad \hbox{when}\quad x\to -\infty,\\
u_{2}(x,k)\sim& \frac{1}{T_{1}(k)}\ee^{-\ii kx}+
\frac{R_{1}(k)}{ T_{1}(k)}\ee^{\ii kx}
\quad\hbox{when}\quad x\to \infty,\\
\end{aligned}
$$
where
$R_1(k)$, $R_2(k)$, $T_1(k)$ and $T_2(k)$ are complex constants 
determined by $c$ 
and $k$.
The matrix
$$
S(k)=\left( \begin{matrix} T_{1}(k) & R_{2}(k)\\ R_{1}(k) & T_{2}(k)\\ \end{matrix} 
\right)
$$ 
is the scattering matrix determined by
$c$; $R_{1}(k)$ and  $R_{2}(k)$ are the reflection coefficients, whereas
 $T_{1}(k)$ and  $T_{2}(k)$ are the transmission coefficients.

We shall assume that  $R_{2}(k)=R_{2}(k;c)$ is known.
Our aim here is to prove the uniqueness part of the inverse scattering 
problem, that is, to prove the next theorem.

\begin{theorem}\label{thm:4.1.1}
The mapping
$$
\begin{gathered}
{\mathcal R}\colon\{c\colon\R\to\R, \text{ {\rm  measurable} }\mid c \ \text{\rm satisfies  [H1] 
and [H2]}\}\rightarrow L^\infty(\R); \\
{\mathcal R}(c)=R_2(k;c) 
\end{gathered}
$$
is injective. 
\end{theorem}

Even in the simple case of \eqref{eq:4.1.1}, the methods used in the Schr\"odinger case
are not available in this case, since 
the behaviour in $k$, $|k|$ large,  of the solutions 
$u_{1}(x,k)$, $u_{2}(x,k)$  is no longer easy to control, 
at least when $c$ is not assumed to be 
smooth enough.
(See \cite{F}, \cite{DT}, \cite{M},  
\cite{AKM}, and the references therein for the methods usd in the Schr\"odinger case,  and also \cite{B01} 
for the same problem in higher dimensions, where the extra dimensions and extra decay allow the reduction to a 
Schr\"odinger type case.) 

We however adapt the idea in \cite[Sect.5, Lemma~1]{DT} for the proof of Levinson's theorem, 
using properties of Nevanlinna functions instead of Hardy classes. 
The main tool in our proof is the uniqueness result contained in Proposition~\ref{uniprop}, which says that under certain conditions, a Poisson type representation holds for
Nevanlinna functions on the upper-half space. 
Namely, it is well known (\cite{AD}) that when $F$ is holomorphic in the upper half-plane and has non-negative real part, 
then $F$ can be written as
\begin{equation}\label{int:1}
F(z) =\alpha z+\beta + \frac{1}{\pi\ii} \int \limits_{\R} \big( \frac{1}{t-z}-\frac{t}{t^2+1}\big) \dd \mu(t), \quad \im z\ge 0, 
\end{equation}
where $\alpha$, $\beta\in \C$, and $\mu$ is a positive measure on $\R$. 
We show that a representation of the type of \eqref{int:1} holds for certain Nevanlinna functions as well, without necessarily knowing 
\textit{a priori} that they have nonnegative real part. 

The paper is organized as follows: Section~\ref{sect:nevanlinna} recalls the necessary notions and notations concerning Nevanlinna functions, 
and gives the Poisson-type representation for certain classes of functions. 
Section~\ref{sect2} recalls the construction of the Jost solutions $u_1(x,k)$, $u_2(x,k)$
briefly presents the scattering matrix, and gives  some of their properties.
Section~\ref{sect:w}  deals with properties of the function $r(x, k)$ (see \cite{SyWG}),  
which is $\ee^{-2\ii kx}$ times the reflection coefficient for the problem \eqref{eq:4.1.1} with $c$ replaced by the 
function $c_x(y) = c(y)$ when $y> x$ and $c_x(y)=1$ otherwise. 
The fact that
$T_1(k)$ is uniquely determined by
$R_2(k)$ is proved in Section~\ref{sect:T}. 
Some further analysis of the behaviour in $k$ of 
$u_{1}(x,k)$ is given in Section~\ref{sect:further}.
The main result (Theorem~\ref{thm:4.1.1}) is proved in  Section~\ref{sect:uniq}, using the 
Poisson-type representation obtained in Section~\ref{sect:nevanlinna}.

\subsection*{Notation} 
We denote 
$[f,g]_x=f'(x)g(x)-f(x)g'(x). $
Note that $-[f, g]_x$ is the Wronskian of  $f$ and  $g$.
If this quantity is constant on  $ x\in \R$, as it happens when 
$f$, $g$ are solutions 
to the equation \eqref{eq:4.1.1}, it will be denoted  by 
$[f,g]$.

Throughout the paper we use the notation $\Pi$ for  the upper half-plane in $\C$, 
$$\Pi= \{z\in \C\mid \im z>0\}.$$
%Also $H^p(\Pi)$, $1\le p\le \infty$ is the Hardy space on $\Pi$,  $N(\Pi)$ denotes the space of 
%Nevalinna functions on $\Pi$ (functions $f\colon \Pi \to \C$ that are analytic and $\log|f(z)|$ has an harmonic majorant), $N^+(\Pi)$ is the subspace  of $N(\Pi)$ consisting of functions of the form $h/g$ where $h$, $g\in H^{\infty}(\Pi)$, $\abs{h}, \abs{g}\le 1$,  and $g$ is outer. 
We refer to \cite{RR} for notation and results on Hardy and Nevanlinna classes.
 Also,  we use the shorthand notation $\langle x\rangle= (1+\vert x\vert^2)^{1/2}$, $x\in \R$.

\section{Nevanlinna functions and a uniqueness lemma}\label{sect:nevanlinna}

\begin{definition}\label{Nevanlinna}

\begin{itemize}
\item[(a)] The Nevanlinna class 
$N(\Pi)$ is the class of those holomorphic functions $F$ defined on $\Pi$ such that $\log^+\abs{F}$ has a harmonic majorant on $\Pi$.

\item[(b)] Recall that $F\in N^+(\Pi)$ if and only if $F=AG$ where $A\in N(\Pi)$ is inner and $G\in N(\Pi)$ is outer.
\end{itemize}
\end{definition}

\begin{remark}\label{rem-nevanlinna}
We recall the following properties of the Nevanlinna functions.
\begin{enumerate} 
\item
 If $F, G\in N(\Pi)$ then $F+G$, $FG\in N(\Pi)$. 
Moreover if $G\not\equiv 0$ and $F/G$ is holomorphic on $\Pi$, then $F/G\in N(\Pi)$. 

\item 
If $F, G\in N^+(\Pi)$ then $F+G, FG\in N^+(\Pi)$.

\item
It can be shown (see \cite[Thm. 3.20]{RR}) that $N(\Pi)$ is the set of quotients $F=G/H$, where $G$ and $H$ are holomorphic and bounded on $\Pi$, and $H$ is not vanishing on $\Pi$.

\item  $N^+(\Pi)$ is the set of quotients $F=H/V$, where $H$ and $V$ are holomorphic and bounded  by $1$ on $\Pi$, and $V$ is outer.

\item $N^+(\Pi)$ is the smallest algebra of functions containing all inner and outer functions in $N(\Pi)$. 
\end{enumerate}
\end{remark}

\begin{remark}\label{im-positive}
If $F\not\equiv 0$ is holomorphic on $\Pi$ and $\im F\ge 0 $ on $\Pi$, then $F\in N(\Pi)$ and it is outer. 
\end{remark}

Every function $F\not\equiv 0$ in $N(\Pi)$ has a factorization 
\begin{equation}\label{factorization}
F(z) =\ee^{-\ii\tau z} B(z) G(z) \frac{S_+(z)}{S_-(z)}, \quad \im z>0,
\end{equation}
where $\tau$ is a real number, $B$ is a Blaschke product, $G$ is an outer function and $S_+$ and $S_-$ are singular inner functions. 
The factorization is essentially unique: $\ee^{-\ii\tau z}$, $S_+$ and $S_-$ are uniquely determined, while $B$ and $G$ are determined up to  multiplicative constants of modulus $1$.

\begin{definition}\label{exptype}
\begin{itemize}
 \item[(a)]
The real number $\tau$ in \eqref{factorization} is called \textit{the mean type} of $F$.
\item[(b)] 
An entire function $F$ is of \textit{exponential type} if 
$$ \tau_F =\limsup_{\abs{z}\to \infty} \frac{\log{\abs{F(z)}}}{\abs{z}} <\infty. $$
The number $\tau_F$ is called the exact type of $F$. 
\end{itemize}
\end{definition}

\begin{remark}\label{kreintype}
We collect here important properties of the mean type and Nevalinna functions
(see \cite[Thm. 6.15--6.17]{RR}).

 \begin{enumerate}
\item\label{meantype} 

Let $0\not\equiv F\in N(\Pi)$. 
Then its mean type $\tau$ is given by
$$ \tau=\limsup_{y\to \infty} \frac{1}{y}\log\abs{F(\ii y)}.$$

\item\label{Krein}
Let $F$ be an entire function. The following assertions are equivalent.
\begin{itemize}
\item[(a)] $F$ is of exponential type and 
$$ \int\limits_{-\infty}^{\infty} \frac{\log^+ \abs{F(t)}}{1+t^2}\, \dd t< \infty.$$
\item[(b)] The restrictions of $F(z)$ and $\widetilde{F}(z):=\overline{F(\bar{z})}$ to the upper half-plane belong to $N(\Pi)$. 
\end{itemize}

\item\label{meantype2}
Let $F\not\equiv 0$ be a function that satisfies the equivalent conditions (a) and (b) in \eqref{Krein}.
Let $\tau_F$ be the exact type of $F$, and let $\tau_+$ and $\tau_-$ be the mean types of the restrictions of $F$ and $\widetilde{F}$ to $\Pi$, respectively. 
Then $\tau_+ + \tau_-$ is nonnegative, and 
$$ \tau_F =\max(\tau_+, \tau_-).$$
\end{enumerate}
\end{remark}

\subsection{A Poisson-type representation}

\begin{lemma}\label{unilemma1}
Let $h \in L^1(\R, \dd t/(1+t^2))$ be a real function. 
Define
$$ f(z) =\frac{1}{\pi \ii} \int h(t) \Big( \frac{1}{t-z} -\frac{t}{t^2+1}\Big)\, \dd t, \quad z\in \Pi.$$
Then $f\in N(\Pi)$.
\end{lemma}

\begin{proof} 
It is clear that $f$ is analytic on $\Pi$. 
Set 
$$ f_+(z) = \frac{1}{\pi \ii} \int \abs{h(t)} \Big( \frac{1}{t-z} -\frac{t}{t^2+1}\Big)\, \dd t, \quad z\in \Pi.
$$
Then $f_+$ is analytic on $\Pi$, and $\re f_+ \ge 0$ on $\Pi$, hence $f\in N(\Pi)$. 
Also $f_+ -f$ is analytic on $\Pi$ and it is easily seen that, when $z=\lambda + \ii \kappa\in \Pi$,  
$$ \re (f_+ -f)(z) =\frac{\kappa}{\pi} \int \frac{\abs{h(t)} -h(t)}{(t-\lambda)^2 + \kappa^2} \, \dd t\ge 0, \quad 
 $$
hence $f_+ -f\in N(\Pi)$.
Thus $f = f_+ -(f_{+} - f)\in N(\Pi)$.
\end{proof}

The next proposition is a Poisson-type representation for a particular class of Nevanlinna function.

\begin{proposition}\label{uniprop}
Let $h$ be a function in $N(\Pi)$, continuous on $\overline{\Pi}$.   
Assume that $\re h/(1+t^2)\in L^1(\R)$, $\re h$ is even on $\R$ and  $h$ is bounded on $\ii \R$.
Then there exist $\alpha$, $\beta\in \C$ such that 
$$ h(z) =\alpha z +  \beta+  
\frac{1}{\pi \ii} \int \re h(t) \Big( \frac{1}{t-z} -\frac{t}{t^2+1}\Big)\, \dd t
$$ when $z\in \Pi$.
 \end{proposition}

\begin{proof}
First note that the boundedness of $h$ on $\ii \R$ implies that $h$ has nonpositive mean type.
Consider the function defined by   
$$ g(z)=\frac{1}{\pi \ii} \int \re h(t) \big( \frac{1}{t-z} -\frac{t}{t^2+1}\big)\, \dd t, \quad z\in \Pi. $$
 Then it is enough to show that $h-g$ is a polynomial of degree $1$.

It follows by the previous lemma that $g\in N(\Pi)$.
Its real part   $\re g$ extends continuously to $\overline{\Pi}$. 
Note that 
$$\im g(\ii \kappa)= \frac{1}{\pi} \int  \re h(t) \Big( \frac{t}{t^2+\kappa^2} +\frac{t}{t^2+1}\Big)\, \dd t=0, \quad  \kappa> 0,$$
since $\re h$ is even. 
It follows that $g(\ii\kappa)/\kappa\to 0$ when $\kappa\to \infty$. 
In particular,  $g$ has nonpositive mean type.

Let us  set $f(z) = \ii [h(z)- g(z)]$. 
Then 
$f $ belongs to $N(\Pi)$, and it has nonpositive mean type.  
Also, $\im f(z)$ tends to $0$ when $\Pi \ni z\to x\in \R$. 
Hence, by the reflection principle (see \cite[Thm. 11.11 and Thm. 11.17]{Ru}), the function
$$ F(z) =\begin{cases} f(z), & z\in \Pi\\
\overline{f(\bar z)},  & \im z < 0,\end{cases} $$  
extends analytically on $\C$. 
It is clear that $F\vert_{\Pi}\in N(\Pi)$, and $\tilde{F}(z) = \overline{F(\bar{z})}= f(z)$, hence $\tilde F\in N(\Pi)$. 
It follows by Remark~\ref{kreintype} that $F$ is of exponential type $0$. 
Since $F/(\ii \kappa+1)$ is bounded on $\ii \R$, it follows that $F$ is a polynomial of degree $1$.
\end{proof}

The next corollary to Proposition~\ref{uniprop} will be used in proving the 
uniqueness in our scattering problem, in Section~\ref{sect:uniq}.

\begin{corollary}\label{cor1001}
Let $h$ be a function in $N(\Pi)$, continuous on $\overline{\Pi}$. 
Assume that
\begin{enumerate}
\item\label{cor1001_1} $h(0)=0$, $\overline{h (z)}= h(-\overline{z})$ when $z\in \overline{\Pi}$,
\item\label{cor1001_2} $h$ is bounded on $\ii \R$, 
\item\label{cor1001_3}
$\re{h}/(t(1+\vert t\vert) ) \in L^1(\R)$. 
%\item\label{cor1001_4} There is $0\le \gamma< 1$ such that $\re{h}/(t^{2}+1)^{\gamma/2}\in L^1(\R)$.
\end{enumerate}
Then
\begin{equation}\label{cor1001_5}
h(z) =
\frac{1}{\pi \ii} \int \re h(t) \Big( \frac{1}{t-z} -\frac{t}{t^2+1}\Big)\, \dd t,  \qquad z\in \Pi.
\end{equation}
for all  $z\in \Pi$,
 \end{corollary}

\begin{proof} 
By Proposition~\ref{uniprop} there are $\alpha$, $\beta\in \C$ such that
\begin{equation}\label{101}
h(z) =\alpha z + \beta + 
\frac{1}{\pi \ii} \int \re h(t) \Big( \frac{1}{t-z} -\frac{t}{t^2+1}\Big)\, \dd t
\end{equation} when $z\in \Pi$.
Since $\re{h}$ is even on $\R$ we get that 
\begin{equation}\label{102}
h(\ii \kappa) = \ii \alpha \kappa + \beta + \frac{\kappa}{\pi}\int\frac{ \re h(t)}{t^2+\kappa^2} \dd t, \quad \kappa>0.
\end{equation}
The hypothesis \eqref{cor1001_3} ensures that 
$$ \int \frac{\re h(t)}{t^2+\kappa^2} \dd t \to 0 \quad \text{when} \; \kappa\to \infty. $$
Hence we see from \eqref{102} that $\alpha=0$. 

On the other hand $h(\ii \kappa)\to 0$ when $\kappa\to 0$, and 
$$ \kappa \int \frac{\re{h(t)}}{t^2+\kappa^2} \dd t \to 0 \quad \text{when} \; \kappa\to 0. $$
Taking  reals parts in \eqref{101} when  $k =\ii \kappa$,  $\kappa>0$,  and then letting $\kappa$ tend to $0$,  we obtain that $\beta=0$. 
\end{proof}

 %%%%%%%%%%%%%%%%%%%%%%%%%%%%%%%%%%%%%%%%%%%%%%%%%%%%
\section{Jost solutions and the scattering matrix}\label{sect2}

We start here the study of the scattering problem for \eqref{eq:4.1.1}, under the conditions [H1] and [H2].
First we recall well-known properties of the Jost solutions and scattering matrix for \eqref{eq:4.1.1}, and 
prove some 
extra simple properties. 

\subsection {Jost solutions}
The scattering matrix 
can  be given in terms of Jost solutions, that is, solutions $(u_{j}(\cdot,k))_{j=1,2}$, $k\in \overline \Pi$, 
of the equation \eqref{eq:4.1.1},  
that satisfy
$u_1(x,k)\sim \ee^{\ii kx}$ when  $x\to \infty$, and  $u_2(x,k)\sim \ee^{-\ii kx}$ when 
$x\to -\infty$.

The results of this subsection are quite  standard and can be proved  as in the
Schr\"odinger case (see \cite{DT}) for the potential $k^2q$, where
\begin{equation}\label{eq:defq}
q=1-\frac{1}{c^2},
\end{equation}
which satifies (from [H1],  [H2])  
$$
q\in L^1(\R) \cap L^\infty(\R).
%|q(x)|\le C\langle x\rangle^{-1-\delta}, \quad a.e. \; x\in \R.
$$
%Set
%$
%\gamma_0=\sup_{x\in \R} |q(x)|$.

\begin{theorem}\label{thm:2.1}
For every $k\in \overline{\Pi}$ there exist unique solutions 
$u_1(x,k)$, $u_2(x,k)$ to the equation
\eqref{eq:4.1.1} such that
\begin{equation}\label{eq:2.2}
\ee^{-\ii kx} u_1(x,k)\to 1 ,\quad \ee^{\ii kx} u_2(x,k) \to 1,
\end{equation}
when $x\to \infty$ and 
$x\to -\infty$, respectively.
The following estimates hold:
\begin{align}
&|u_1(x,k)-\ee^{\ii kx}|\le |k|\gamma(x)\ee^{|k|\gamma(x)-x\im k}, \label{eq:2.4i}\\
&|u_2(x,k)-\ee^{-\ii kx}|\le |k|\eta(x)\ee^{|k|\eta(x)+x\im k}, \label{eq:2.4ii}\\
& |u'_1(x,k)-\ii k u_1(x,k)|\le |k|^2\gamma(x)\left(1+|k|\| q\|_{L^1}\ee^{|k|\| q\|_{L^1}}\right)
\ee^{-x\im k},\label{eq:2.5i}\\
& |u'_2(x,k)+\ii k u_2(x,k)|\le |k|^2\eta(x)\left(1+|k|\| q\|_{L^1}\ee^{|k|\| q\|_{L^1}}\right)
\ee^{x\im k}, \label{eq:2.5ii}\\
&\overline{u_1(x,k)}=u_1(x,-\overline{k}),\qquad \overline{u_2(x,k)}=u_2(x,-\overline{k}), 
\label{eq:2.6}
\end{align}
where $\gamma(x)=\int_{x}^{\infty}|q(y)|\dd y$ and 
 $\eta(x)=\int^{x}_{\infty}|q(y)|\dd y$.

When $x\in \R$ is fixed, the functions
$$
\Pi\ni k\mapsto u_j(x,k)\in \C, \qquad \Pi \ni k\mapsto u'_j(x,k)\in \C,
\quad j=1, 2,$$ 
are holomorphic and extend continuously at $\im k=0$.
\end{theorem}
%%%%%%%%%%%%%%%%%%
We denote 
\begin{equation}\label{def:m}
 m_1(x, k) = \ee^{-\ii kx} u_1(x, k), \quad 
m_2(x, k) = \ee^{\ii kx} u_2(x, k).
\end{equation}

\begin{remark} \label{rem2.2}

The functions 
$m_1(x, k)$ and $m_2(x, k)$  solve the equations
\begin{align}
m_1''(x, k)+ 2\ii k m_1'(x, k) & = k^2 q(x) m_1(x, k), \label{eq:m1} \\ 
m_2''(x, k)- 2\ii k m_2'(x, k) & = k^2 q(x) m_2(x, k),\label{eq:m2}
\end{align}
respectively.
Theorem~\ref{thm:2.1} gives  
$$
\begin{gathered}
|m_1(x, k) -1|\le |k|\gamma(x)\ee^{|k|\gamma(x)},\quad 
|m_1'(x, k)| \le |k|^2\gamma(x)\left(1+|k|\| q\|_{L^1}\ee^{|k|\| q\|_{L^1}}\right),\\
 |m_2(x, k) -1|\le |k|\eta(x)\ee^{|k|\eta(x)}, \quad
 |m_2'(x,k)|\le |k|^2\eta(x)\left(1+|k|\| q\|_{L^1}\ee^{|k|\| q\|_{L^1}}\right)
\end{gathered}
$$ 
when $k\in \overline{\Pi}$ and $x\in \R$.
Note that 
\begin{equation}\label{m1:ser}
m_1(x, k)   = 1+ \sum\limits_{n=1}^\infty g_n(x, k)
\end{equation}
with
$$
g_n(x, k) = \Big(\frac{k}{2\ii}\Big)^n \idotsint\limits_{x\le x_1\le \dots \le x_n} 
D_k(x_1-x) \cdots D(x_n-x_{n-1})  q(x_1) \cdots q(x_n) \dd x_1 \dots \dd x_n, \nonumber
$$
where $D_k(x) = \ee^{2\ii kx}-1$. 
\end{remark}

\subsection{The scattering matrix}

Let $u_{1}(\cdot,k)$, $u_{2}(\cdot,k)$ 
be the solutions introduced in the previous subsection.
Consider $k\in \R\setminus\{0\}$.
It follows from  \eqref{eq:2.4i}, \eqref{eq:2.5i}  
that
\begin{equation}\label{eq:2.7}
\comm{u_1(\cdot,k)}{u_1(\cdot,-k)} 
 =2\ii k\ne 0, \\
\end{equation}
and
similarly,  
\begin{equation}\label{eq:2.9}
\comm{u_2(\cdot,k)}{  u_2(\cdot,-k)}=-2\ii k\ne 0, 
\end{equation}
 when $k\in \R\setminus \{0\}$.
 
We deduce from \eqref{eq:2.6}, that $u_{1}$, $\overline{u}_1$ 
($u_{2}$, $\overline{u}_2$) are linearly independent solutions of
\eqref{eq:4.1.1}. 
Also, it follows that $u_1(\cdot, k)$ and $u_2(\cdot, k)$ are linearly independent solutions 
of \eqref{eq:4.1.1} when $k\ne 0$, due again to $\eqref{eq:2.6}$, \eqref{eq:2.7} and \eqref{eq:2.9}.
Hence there are constants
 $T_{1}(k)$, $T_{2}(k)$,
$R_{1}(k)$, $R_{2}(k)$ such that $1/T_j(k) \ne 0$, $j=1, 2$,  and 
\begin{align}
u_{2}(x,k)&=\frac{R_{1}(k)}{ T_{1}(k)}u_{1}(x,k)+\frac{1}{T_{1}(k)}u_{1}(x,-k),\label{eq:2.11i}\\
u_{1}(x,k)&=\frac{R_{2}(k)}{T_{2}(k)}u_{2}(x,k)+\frac{1}{T_{2}(k)}u_{2}(x,-k), \label{eq:2.11ii}
\end{align}
when $k\ne 0$ is real.   
Here $R_1(k)$, $R_2(k)$, $T_1(k)$, $T_2(k)$ are the scattering coefficients
of  
$c$ at energy $k^2$.
The first two quantities are the reflection coefficients, 
while the others are the transmission coefficients, and
the matrix
$$S(k)=\left(
\begin{matrix} T_1(k) & R_2(k)\cr
                    R_1(k) &T_2(k)\cr
\end{matrix}
\right)$$
is the scattering matrix.
 We see from \eqref{eq:2.7}, \eqref{eq:2.9}, 
\eqref{eq:2.11i} and \eqref{eq:2.11ii} that
\begin{align}
[u_1(\cdot,k), u_2(\cdot,k)]&=\frac{2\ii k}{T_{1}(k)}=\frac{2\ii k}{T_{2}(k)},\label{eq:2.12} \\
[u_2(\cdot,k), u_1(\cdot,-k)]&= 2\ii k\frac{ R_1(k)}{T_{1}(k)}, \label{eq:2.13} \\
[u_2(\cdot,-k), u_1(\cdot,k)]&=2\ii k\frac{R_{2}(k)}{T_{2}(k)}. \label{eq:2.14}
\end{align}
It follows that
\begin{gather}
T_{1}(k)=T_{2}(k){\buildrel{\text{not.}}\over{=\hskip-2pt=}} T(k)\quad \text{and}\quad 
\overline{T(k)}=T(-k),
\label{eq:2.15}
\\
R_{1}(k)T(-k)+R_{2}(-k)T(k)=0,\label{eq:2.16}\\
\overline{R_{1}(k)}=R_{1}(-k),\quad \overline{R_{2}(k)}=R_{2}(-k),\label{eq:2.17}
\\
|T(k)|^{2}+|R_{2}(k)|^{2}=|T(k)|^{2}+|R_{1}(k)|^{2}=1  \label{eq:2.18}
\end{gather}
when $k\in \R$. 

As a consequence, the matrix $S(k)$ is unitary and
\begin{equation}\label{eq:2.20}
|T(k)|\le 1,\quad |R_{1}(k)|\le 1,\quad  |R_{2}(k)|\le 1,
\quad k\in \R\setminus \{0\}.
\end{equation}

The function $2\ii k/T(k)$ can be analytically extended to $\Pi$, taking
\eqref{eq:2.12} as definition.
It extends continuously to $\R\setminus \{0\}$, and  can not vanish at any point
in $\overline{\Pi} \setminus \{ 0\}$. 
Indeed,  on $\R\setminus \{0\}$ this follows from the discussion above, 
while for $k\in \Pi$ this is a consequence of the fact that there are no negative eigenvalues for the operator
$ u\to -c (cu)''$ considered as a self-adjoint operator in 
$L^2(\R)$ with domain $\{ u \in L^2 \mid (cu)'' \in L^2\}$.
We thus see that $T(k)$ can be extended to $\Pi$  as a holomorphic function, 
continuous on $\R^+\setminus \{0\}$. 
In addition $\overline{T(k)}= T(-\overline{k})$ when $k\in {\overline{\Pi}} \setminus 0$.

We have 
(see \cite{DT})
\begin{align}
\frac{R_{2}(k)}{T(k)} &= \frac{k}{ 2\ii }\,\int\limits_{-\infty}^{\infty}
                         \ee^{2\ii kt}\,q(t)\,m_{1}(t,k)\dd t, \label{eq:2.21} \\
\frac{1}{T(k)}&=1-\frac{k}{ 2\ii }\,\int\limits_{-\infty}^{\infty}q(t)\,m_{1}(t,k)\dd t
\label{eq:2.22}
\end{align}
when $k\in \R\setminus\{0\}$.

We need an estimate for $R_2(k)$ when $k$
is in a neighbourhood of $0$.
We have the following lemma.

\begin{lemma}\label{lemma:2.2}
\rm{(i)} The function $\R\setminus \{0\}\ni k\mapsto T(k)\in \C$ extends  
continuously to  $\overline{\Pi}$ with $T(0)=1$.

(ii)  The function
$\R\setminus \{0\}\ni k\mapsto R_2(k)\in \C$ 
extends continuously  to 
 $\R$ with $R_2(0)=0$.
Moreover
\begin{align} 
|R_{2}(k)|& \le C|k| \quad\text{when}\quad k\in \R, \label{eq:2.23}\\ 
\lim_{k\to 0} \frac{R_{2}(k)}{k} & =\frac{1}{2\ii}\,
\int\limits_{-\infty}^{+\infty}q(t)\dd t ,\label{eq:2.24}
\end{align}
where $C$ is a constant independent on $k$.
\end{lemma}
\begin{proof}
The first part of the proof follows from the discussion before the present lemma, relation \eqref{eq:2.22}, 
and Theorem~\ref{thm:2.1}.

It follows from \eqref{eq:2.21} that $k\to R_2(k)/T(k)$ extends continuously on $\R$,
and since  $k\to T(k)$ extends continuously  on 
$\R$
we get that $k\to R_2(k)$ has the same property.
The remaining assertions of the statement are direct consequences of 
\eqref{eq:2.21} and \eqref{eq:2.22}, 
taking into account that
$|T(k)|\le 1$  and
\eqref{eq:2.4i} holds.
\end{proof}

The following simple lemma will be essentially used in the paper.

\begin{lemma}\label{slemma}
Let $h_1$, $h_2$, $g_1$, $g_2\in \C$, $a$, $b\in \C$ with $\abs{a}^2 +\abs{b}^2=1$. 
Assume that 
\begin{align}
a h_1  & = bh_2 + \bar{h}_2\label{1001},\\
a g_1 & =bg_2 + \bar{g}_2\label{1001a}
\end{align}
Then 
$
2\re (ah_1 g_2)= \abs{a}^2 h_1\bar{g}_1 + \abs{a}^2h_2\bar{g}_2
$.
\end{lemma}

\begin{proof} 
We may assume $g_2\ne 0$ since otherwise the equality to prove obviously holds.
Multiplying \eqref{1001} by $g_2$ and taking real parts,  we get
\begin{equation}\label{real}
2\re(a h_1 g_2) =2\re (bh_2 g_2) + \bar{h}_2g_2+  h_2\bar{g}_2.
\end{equation}
On the other hand by multiplying \eqref{1001} and the conjugate of \eqref{1001a} it follows that
$$ 
\begin{aligned}
\abs{a}^2h_1\bar{g}_1 & = (bh_2 + \overline{h}_2)(\bar{b} \bar{g}_2 + g_2)\\
& =\abs{b}^2 h_2\bar{g}_2 + 2\re (bh_2g_2) + \bar{h}_2 g_2\\
& = (1-\abs{a}^2) h_2\bar{g}_2 + 2\re (bh_2g_2) + \bar{h}_2 g_2\\
&= -\abs{a}^2 h_2\bar{g}_2+ 2\re (ah_1g_2),
\end{aligned}
$$
where we have used  the equality $\abs{a}^2 +\abs{b}^2=1$ and \eqref{real} to express $2\re (bh_2g_2)$.
\end{proof}

We use Lemma~\ref{slemma} for $a=T(k)$, $b= R_2(k)\ee^{-2\ii k x}$, $h_1= g_1=m_1(x, k)$ and
 $h_2= g_2=m_2(x, k)$ to get the next corollary.
\begin{corollary}\label{lemma:2.3}
The equality
\begin{equation}\label{local:201}
2\re (T(k) m_1(x, k)m_2(x, k))
=|T(k)m_1(x, k)|^2+ |T(k) m_2(x, k)|^2.
\end{equation}
holds for every  $x\in \R$ and $k\in \R$.
\end{corollary}
%
%

%%%%%%%%%%%%%%%%%%%%%%%%%%%%%%%%%%%%%%%%%%%%%%%%%%%%%%%%%%%%%%%%%%%

\section{The functions $w$ and $r$}\label{sect:w}

We consider the functions defined by
\begin{equation}\label{eq:defw}
w(x,k)=\begin{cases}
\frac{u_1'(x,k)}{\ii k u_1(x,k)} & \text{when}\quad   k\in \overline{\Pi} \setminus \{0\}\cr
                 1 & \text{when}\quad k=0\cr
\end{cases},\quad x\in \R,
\end{equation}
and
\begin{equation}\label{eq:3.1}
r(x,k)=\frac{1-w(x,k)}{1+w(x,k)}, \quad x\in\R, \; k\in \overline{\Pi}, 
\end{equation}
Note that $w(x, k)$ is a version of a Titchmarsh-Weyl function for the problem \eqref{eq:4.1.1} on $[x, \infty)$ with $c$ replaced by 
$c_x(y)= c(y)$, $y\ge x$, 
and we shall see that $r(x, k)$ is in fact $\ee^{-\ii k x} $ times the reflection coefficient for the same problem. 
We prove here some results for these function, that will be used later on.

We first need to show that the definitions   
above make sense.

\begin{lemma}\label{lemma:3.1}
Let 
$k\in \overline{\Pi}$ be fixed. 
Then $u_1(x,k)\ne 0$ when $x\in \R$.
\end{lemma}

\begin{proof}
Assume  $u_1(x_{0}, k)=0$ for an  $x_{0}\in \R$ ($x_0$ may depend on $k$).
Due to \eqref{eq:2.7} and \eqref{eq:2.6}, this cannot happen 
when  $k\in \R\setminus \{0\}$. 
If $\im k\ne 0$, then $\frac{1}{c}u_1(\cdot, k)$ 
is an eigenfunction  corresponding to the eigenvalue
$k^{2}$ for the operator $H_{D} f=-c(c\,f)''$ when 
$f\in D(H_{D})=\{g\in L^{2}(x_{0},\infty)\mid -c(c\,g)''\in  L^{2}(x_{0},\infty),\,
(cg)(x_{0})=0\}$.
Thus $k^{2}\in\sigma(H_{D})$, which contradicts the fact that $H_{D}\ge 0$.
In view of \eqref{eq:2.4i}, $u_1(x, 0)=1$ for any $x$ real, and this 
concludes the proof of the lemma.
\end{proof}

Lemma~\ref{lemma:3.1}. shows that definition \eqref{eq:defw} 
is correct. 
The function $w(\cdot, k)$   
is locally of class
$W^{1,\infty}(\R)$ (that is, it is continuous and has a locally $L^\infty$ derivative) 
and
satisfies  
\begin{equation}\label{eq:eqw}
w'(x,k)=\frac{\ii k}{c^2(x)}-\ii k w^2(x,k), \quad x\in \R, \; k\in \overline{\Pi}\setminus\{0\}.
\end{equation}
We note that  (by \eqref{eq:2.5i}, \eqref{eq:2.4i}), 
$$
\lim_{k\to 0}\frac{u_1'(x, k)}{ \ii k u_1(x, k)}
=\lim_{k\to 0}\frac{u_1'(x,k)-\ii k u_1(x,k)}{\ii ku_1(x,k)}+1=1.
$$
It follows that 
the function $\Pi\ni k \mapsto w(x,k)\in \C$ 
is analytic and 
continuous on $\overline{\Pi}$  when $x\in \R$ is fixed (see Thm.~\ref{thm:2.1} and Lemma~\ref{lemma:3.1}).

We can also  define
\begin{equation}\label{eq:defw2}
w_- (x,k)=-\frac{u_2'(x,k)}{\ii k u_2(x,k)}\qquad \text{when}\quad x\in\R,\, k\in\overline{\Pi}
\setminus\{0\},
\end{equation}
and show that it has properties similar to those of $w$.

\begin{lemma}\label{lemma:3.2}
The real part of  $w(x, k)$, $x\in \R$, $k\in \overline{\Pi}$ is positive. 
Specifically, 
\begin{equation}
\label{eq:3.4}
\re w(x,k)=\begin{cases} 
            \frac{1}{|u_1(x,k)|^2} & k\in \R, \\ 
\frac{\im k}{|k|^2\,|u_1(x,k)|^2}
\int\limits_x^{\infty}\left(\frac{|k|^2 |u_1(y,k)|^2}{c^2(y)}+|u_1'(y,k)|^2\right)\,\dd y, 
&  k\in \Pi, 
\end{cases}
\end{equation}
for every $x\in \R$.
\end{lemma}

\begin{proof}
If  $k\in \R\setminus \{0\}$, the equality \eqref{eq:3.4} is a consequence of
\eqref{eq:2.6} and
\eqref{eq:2.7}. 
When $k=0$ both sides in \eqref{eq:3.4} are equal to $1$.
\par
 Assume now that $\im k>0 $. We have
\begin{equation}\label{eq:3.6}
\re w=\frac{1}{2\ii |k|^2\,|u_1|^2}(u_1'\,\overline{u_1\,k}- 
\overline{u_1'}\,u_1 \,k).
\end{equation}
On the other hand, since  $u_1$ solves \eqref{eq:4.1.1} and $u_1$ and $u'_1$ decay
exponentially when $ x\to \infty$, we may write
$$
\begin{aligned}
(\overline{u_1} u_1')(x,k)&=-\int\limits_x^\infty \overline{u_1}(y,k)\, u_1''(y,k)\dd y-\int
\limits_x^\infty |u_1'(y,k)|^2\dd y\\
&=\int\limits_x^\infty \left(\frac{k^2|u_1(y,k)|^2}{c^2(y)}-|u_1'(y,k)|^2\right)\dd y.
\end{aligned}
$$
Using this equality in \eqref{eq:3.6} we get \eqref{eq:3.4} in this case.
\end{proof}

\begin{remark}
Similarly one may prove  that 
$$\re w_-(x, k)> 0 \quad\text{when} \quad \im k\ge 0.$$
\end{remark}

\begin{corollary}\label{cor:poz}
Let $x\in \R$ be arbitrarily fixed. Then
the function
\begin{equation}\label{corpoz:1} 
k  \longrightarrow \frac{1}{k^2+1}( \vert T(k) m_1(x, k)\vert ^2 
+ \vert T(k) m_2(x, k)\vert ^2 )  
\end{equation}
is integrable on $\R$.
\end{corollary}

\begin{proof}
Note  that by \eqref{eq:2.12} we have 
$$ \frac{2}{T(k)m_1(x, k)m_2(x, k)}= \frac{2}{T(k)u_1(x, k)u_2(x, k)}= w(x, k) + w_{-}(x, k), \quad k \in \overline{\Pi}.
$$
Then by  \eqref{local:201} it follows that  
that the function  in \eqref{corpoz:1} is  $(1+k^2)^{-1}$ times   
the real part of 
$$
\R\ni k \rightarrow   \frac{1}{w(x, k)+w_-(x, k)}\in \C.
$$
The statement follows from the fact that this function, which is holomorphic on $\Pi$,  has positive real part and  extends 
continuously to $\im k=0$.
(See [AD].)
\end{proof}

In the next proposition we use the notation $\gamma_0=\opn{ess sup}_{x\in \R} \abs{q(x)}$.

\begin{proposition}\label{prop:3.3}
When $x\in \R$ and $\ka>0$ one has
\begin{equation} \label{eq:3.7}
\frac{2}{2+c_M^2\gamma_0}\le w(x, \ii \ka)\le 1+\frac{1}{2}\gamma_0.
\end{equation}
\end{proposition}

\begin{remark}
Note that $w(x, \ii \ka)$ is real when $x\in \R$ and $
\ka\ge 0$, since
$\overline{u_1(x, \ii \ka)}=u_1(x, \ii\ka)$.
\end{remark}

\begin{proof}
We set $v=v(\cdot, \ii\ka)=w(\cdot, \ii\ka)-1$.
We get then from \eqref{eq:eqw} and  
\eqref{eq:2.4i} that
$$
\begin{gathered}
v(x, \ii \ka)\to 0 \quad \text{when} \quad x\to +\infty,\\
v'=\ka q+\ka v^2+ 2\ka v\ge \ka q+2\ka v.
\end{gathered}
$$
It follows that
\begin{equation}\label{eq:3.8}
v(x, \ii \ka)\le -\int\limits_x^\infty \ka q(y)\, \ee^{-2\ka (y-x)}\, \dd y\le 
\frac{1}{2}\sup_{t\in [x, \infty)}|q(t)|\le \frac{1}{2}\gamma_0,
\end{equation}
hence the second inequality in \eqref{eq:3.7} holds.

On the other hand 
$\tilde{w}= \tilde{w}(x, \ii\ka)=(w(x, \ii \ka))^{-1}$
($w(x, \ii \ka)>0$)  satisfies 
$$
\tilde{w}'=\frac{\ka}{c^2} \tilde{w}^2-\ka, 
$$
and $\tilde{w}(x, \ii\kappa) \to 1$ when $x\to \infty$.
We set $\tilde{v}=\tilde{w}-1$ and get
$$
\tilde{v}'=
 \frac{\ka}{c^2} \tilde{v}^2+\frac{2\ka}{c^2}\tilde{v}-\ka q\ge 
\frac{2\ka}{c^2}\tilde{v}-\ka q. 
$$
Hence
$$
\tilde{v}(x, \ii \ka)\le  \int\limits_x^\infty \ka q(y)\,
\ee^{-2\ka\int\limits_x^y \frac{1}{c^2(s)}\dd s}\, \dd y
\le \int\limits_x^\infty \ka |q(y)|\,
\ee^{-2\ka (y-x) /c_M^2}\, \dd y 
\le \frac{1}{2}c_M^2\gamma_0.
$$
Thus we have obtained 
\begin{equation}\label{eq:3.10}
 \frac{1}{w}\le 1+ \frac{1}{2}c_M^2\gamma_0,
\end{equation}
which proves the first inequality in \eqref{eq:3.7}.
\end{proof}

\begin{lemma}\label{lemma:3.4}
Let  $c_1$ and $c_2$ be two real measurable functions satisfying 
[H1] and [H2], and let  $q_1$, $q_2$ and 
$w_1$, $w_2$ be the functions  defined  by 
\eqref{eq:defq}, \eqref{eq:defw}, corresponding to $c_1$, $c_2$, respectively.
Then
$$
\| (w_1-w_2)(\cdot, \ii \ka)\|_{L^1(\R)}\le \frac{1}{\alpha}
\|q_1-q_2\|_{L^1(\R)} \quad\text{when} \quad \ka\ge 0,
$$
where
$$
\alpha=\frac{2}{2+c_{M,1}^2\gamma_{0, 1}}+ \frac{2}{2+c_{M,2}^2\gamma_{0,
2}}.
$$ 
Here, for $j=1$, $2$,  $c_{M, j}$ is the constant in [H1] corresponding to
$c_j$, while $\gamma_{0, j}=\sup_t |q_j(t)|$.
\end{lemma}

\begin{proof}
We see that  $w_1(\cdot, \ii \ka)-w_2(\cdot, \ii \ka)$ satisfies
$$
\begin{aligned}
&(w_1(\cdot, \ii \ka)-w_2(\cdot, \ii \ka))' =\ka (q_1-q_2)+\ka (w_1(\cdot, \ii
\ka)-w_2(\cdot, \ii \ka))
(w_1(\cdot, \ii \ka)+w_2(\cdot, \ii \ka)),\\
& w_1(x, \ii \ka)-w_2(x, \ii \ka)\to 0\quad \text{when}\quad x\to \infty.
\end{aligned}
$$
Hence we get 
$$ 
(w_1-w_2)(x, \ii \ka)=\ka \int\limits_x^\infty
(q_2-q_1)(y)\ee^{-\ka\int\limits_x^y(w_1+w_2)(s, \ii \ka)\dd s}
\dd y.
$$
Proposition~\ref{prop:3.3} ensures that
$(w_1+w_2)(s, \ii \ka)\ge \alpha$,
hence
$$
|(w_1-w_2)(x, \ii \ka)|\le \ka \int\limits_x^\infty
|(q_1-q_2)(y)|\ee^{-\alpha\ka(y-x)}
\dd y=\ka \int\limits_0^\infty
|(q_1-q_2)(x+y)|\ee^{-\alpha\ka y} \dd y.
$$
The lemma follows  by integrating this inequality with respect to $x$.
\end{proof}

Lemma~\ref{lemma:3.2} shows that
 $\re (1+w(x,k))>1$, 
therefore $r(x, k)$ may be defined by \eqref{eq:3.1}.
It is easily seen that when $x\in \R$ is arbitrarily fixed, 
the function $\Pi\ni k\mapsto r(x,k)\in\C$
is analytic  and extends continuously at $\im k=0$.

The basic properties of  $r$ are contained in the next lemma. 
%%%%%%%%%%%
\begin{lemma}\label{prop:3.3b}
The function $r(\cdot,k)$ defined by  \eqref{eq:3.1} 
satisfies:
\begin{itemize}
\item[\rm (i)]
 $|r(x,k)|<1$ when $x\in \R$, $k\in \overline{\Pi}$.
\item[\rm (ii)] When $k\in \overline{\Pi}$ is fixed 
%$r(\cdot,k)$ belongs to the Sobolev space $W^{1,\infty}(\R)$ and 
\begin{align}
&r'(\cdot, k) =-2\ii k\,r(\cdot, k) +\frac{\ii kq}{2} (1+r(x, \cdot, k))^{2},\label{eq:3.7b}\\
&\lim_{x\to\infty}r(x,k)=0\quad \text{when}\quad k\in\overline{\Pi},\label{eq:3.8b}
\end{align}
In addition,
\begin{equation}\label{eq:3.9b}
\lim_{x\to -\infty}\ee^{2\ii kx}r(x,k)=R_2(k)\quad \text{when}\quad k\in\R.
\end{equation}
\item[\rm (iii)] If $x\in \R$ is arbitrarily fixed, the function $\overline{\Pi}\setminus \{0\} \ni k\mapsto r(x, k) /k \in \C$ 
extends continuously to $\overline{\Pi}$.
\end{itemize}
\end{lemma}
%%%%%%%%%%%%%%%%%%%%%%%%
\begin{proof}
The inequality in (i) is straightforward  since $r= (1-w)/(1+w)$ and 
$\re w>0$ (Lemma~\ref{lemma:3.2}).

We have 
\begin{equation}\label{eq:3.10b}
r'=-\frac{2}{(1+w)^2}w',\
\end{equation}
$w'$ is locally bounded, $|1+w|>1$, hence $r'(\cdot, k)$ is locally bounded.
The equation \eqref{eq:3.7b} follows from \eqref{eq:3.10b} and  
\eqref{eq:eqw}, while
\eqref{eq:3.8b} is a consequence of  \eqref{eq:2.5i} and \eqref{eq:2.4i}. 

Let $k\in \R\setminus\{0\}$ be fixed. Then by    
\eqref{eq:2.11i}
we have
$$
\ee^{2\ii kx}r(x,k)=
\frac{2\ii k R_2(k)m_2(x,k)-R_2 (k)m'_2(x,k)-\overline{m'_2(x,k)}\ee^{2\ii k x}}
{2\ii k\overline{m_2(x,k)}+R_2 (k)m'_2(x,k)\ee^{-2\ii k x}+\overline{m'_2(x,k)}}.
$$
This gives \eqref{eq:3.9b}, since  
$m_2(x,k)\to 1$ and $m'_2(x,k)\to 0$ when $x\to -\infty$ (see Theorem~\ref{thm:2.1}). 
The equality  \eqref{eq:3.9b} is obvious for $k=0$, since 
$r(x,0)=R_2 (0)=0$ for every $x$.

From \eqref{eq:3.7b} and \eqref{eq:3.8b}  we see that
$$ \frac{r(x, k)}{k} = -\frac{\ii}{2} \int\limits_x^\infty \ee^{2\ii k(y-x)} q(y) (1+ r(y, k))^2 \, \dd y $$ 
when $\im k\ge 0$, $k\ne 0$. Then (iii) follows.
\end{proof}

\begin{remark} 
Let $x\in \R$ be fixed. 
Then $\ee^{2\ii kx} r(x, k)$ is the reflexion coefficient $R_2(k; c_x)$ for $c_x(y)= c(y)$ when $y \geq x$ 
and $c_x(y)=1$ when $y< x$.
\end{remark}

\begin{lemma}\label{lemma:4.5}
We have
\begin{equation}\label{eq:4.9}
\lim_{\kappa\to \infty} qr(\cdot, \ii \kappa)=Q^2
\end{equation}
in $L^1(\R)$.
\end{lemma}

\begin{proof}
(i) Assume first that $c\in C^\infty(\R)$ and obeys 
[H1] and [H2].
We change coordinates
$$ x(y)= \int\limits_0^y\frac{1}{c(s)} \, \dd s $$
and write $v(x(y), k)=u_1(y, k)$, with $u_1$ as in Section~\ref{sect2}.
Then
$$r(y, k)=\Big ( 1-\frac{1}{c(y)}\ \frac{v'(x(y), \ii \ka)}{\ii \ka\ v(x(y), \ii \ka)}\Big)   
\Big ( 1+\frac{1}{c(y)} \ \frac{v'(x(y), \ii \ka)}{\ii \ka\ v(x(y), \ii \ka)}\Big)^{-1}.$$
Since  $c\in C^\infty$, we have that 
$v'(x(y), \ii \ka)/(\ii \ka\ v(x(y), \ii \ka))$ converges 
uniformly  to $1$ when $\ka\to \infty$.
(See \cite{DT}.) 
We get that, for $y$ fixed, 
$$ \lim_{\ka\to \infty} r(y, \ii\ka)=\frac{c(y)-1}{c(y)+1}.$$
It follows that $q(x) r(x, \ii \ka)\to Q^2(x)$ for every  $x$, and 
\eqref{eq:4.9} 
is a consequence of Lebesgue's convergence theorem.

(ii) Assume now that $c$ is a real measurable function that satisfies 
[H1] and [H2].
We set $c_\epsilon=c\ast \varphi_\epsilon= 1+(c-1)\ast \varphi_\epsilon$ and $q_\epsilon=1-c_\epsilon^{-2}$,
 where  $\varphi$ is a nonnegative smooth compactly supported function
with $\int\varphi \, \dd x=1$ and
$\varphi_\epsilon(x)=\frac{1}{\epsilon}\varphi(\frac{x}{\epsilon})$.
Then
\allowdisplaybreaks
\begin{gather}
c_\epsilon(x)\geq c_0\int \varphi_\epsilon(y) \, \dd y =c_0, \nonumber\\
c_\epsilon(x)\leq c_M \int \varphi_\epsilon(y) \, \dd y =c_M, \nonumber\\
|q_\epsilon(x)| \leq \frac{c_\epsilon(x)+1}{c_\epsilon^2(x)} \int |c(x-\epsilon y)-1|\varphi(y) \, \dd y \le \frac{(c_M+1)^2}{c_0^2}=:\gamma_1
\end{gather}
and 
$$ \Vert q_\epsilon-q\Vert_{L^1}\le \frac{1}{c_0^2}\Vert c_\epsilon^2 -c^2\Vert _{L^1}
\to 0 \quad \text{when $\epsilon\to 0$}.
$$
We denote by  $w_\epsilon$ and $r_\epsilon$
the functions  defined as in \eqref{eq:defw} and \eqref{eq:3.1}
corresponding to $c_\epsilon$, and $Q_\ep= 1-c_\epsilon^{-1}$.
Then, since $w_\epsilon(x, \ii \ka) >0$, $w(x, \ii \ka)>0$
and by  lemma~\ref{lemma:3.4}, it follows that there exists  $C$ independent of 
$\epsilon$ and $\ka$ ($C$ may depend on  $c_0$, $c_M$, $\gamma_0$ and $\gamma_1$)
such that
$$
\int\limits|r_\epsilon(x, \ii \ka)-r(x, \ii \ka)|\ \dd x\le 
\int\limits|w_\epsilon(x, \ii \ka)-w(x, \ii \ka)|\ \dd x\le
C\Vert q_\epsilon -q\Vert_{L^1}.
$$
We obtain
$$
\begin{aligned}
\int|q(x)r(x, \ii \ka)-Q^2(x)|\ \dd x  & \le 
\int|q(x)r(x, \ii \ka)-q_\epsilon (x) r_\epsilon(x, \ii \ka)|\ \dd x  \\
+\int & |Q^2_\epsilon(x) -Q^2(x)| \, \dd x+ 
\int|Q^2_\epsilon(x) -q_\epsilon (x) r_\epsilon(x, \ii \ka)| \, \dd x\\
& \le C\Vert c-c_\epsilon\Vert_{L^1} +\int |Q_\epsilon^2(x)-q_\epsilon(x) r_\epsilon(x, \ii \ka)|\ \dd x,
\end{aligned}
$$ 
where $C$ is a constant independent  of  $\kappa$.
Let $\de >0$ be fixed. 
There exists $\epsilon_0=\epsilon_0(\de)$
with $C\Vert c-c_\epsilon\Vert_{L^1}\le \de/2$.
On the other hand, the discussion in (i) shows that there exists $\ka_0=\ka_0(\de)$ such that
if $\ka\geq \ka_0$ then
$$\int |Q_\epsilon^2(x)-q_\epsilon(x) r_\epsilon(x, \ii \ka)|\ \dd x\le \de/2.
$$
Hence, if  $\ka\geq \ka_0=\ka_0(\de)$, we have
$$\int|q(x)r(x, \ii \ka)-Q^2(x)|\ \dd x \le 
\de \ . $$
Since $\de$ has been arbitrarily chosen, this completes the proof of the lemma.
\end{proof}

\begin{theorem}\label{lemma:4.6}
We have
\begin{align}
\int\limits_{\R} \frac{-\log (1-\vert r(x, k)\vert^2 )}{k^2} \dd k  & =\pi \int \limits_x^\infty Q^2(y) \dd y, \quad x\in \R,   
\label{plancherel:1}
\\
\int\limits_{\R} \frac{-\log (1-\vert R_2(k)\vert^2 )}{k^2} \dd k  & =\pi \int \limits_{-\infty}^\infty Q^2(y) \dd y. 
\label{plancherel:2}
\end{align}
In particular, the quantity
$\int_\R Q(x)\ \dd x$
is uniquely determined by the reflection coefficient.
\end{theorem}

\begin{proof}
Let first   $k\in \R\setminus\{ 0\}$ be fixed.
We multiply 
the equation 
$$r'=-2\ii k r+\frac{\ii k q}{2}(1+r)^2,$$
by $\overline{r}$, and  take the real parts.
We obtain
$$(1-|r|^2)'=\frac{i k q}{2}(1-|r|^2)(r-\overline{r}).$$
Since $\lim_{x\to \infty} |r(x, k)|^2=0$ this yields
\begin{equation}\label{eq:4.10}
\frac{-\log (1- |r(x, k)|^2)}{k^2} = 
-\frac{1}{2}\int\limits_x^\infty \frac{q(y)(r(y, k)- \overline{r(y, k)})}{\ii k} \dd y .
\end{equation}
Denote 
$$ h(x, k)= -\int\limits_x^\infty \frac{q(y) r(y, k)}{\ii k} \dd y.$$
For $x$ fixed this is an $H^2(\Pi)$ function with respect to the  variable $k$, with
\begin{equation}\label{eq:4.101}
\overline{h(x, k)}= h(x, -k) \quad \text{when}\quad k\in \R.
\end{equation}
When $x\in \R$, $r(x, k)\ee^{-2\ii k x}$  is the reflection coefficient for $c_x(t) = c(t)$, $t\ge x$, $c_x(t)=1$, $t<x$. 
Hence Proposition~\ref{prop:4.1} and  \eqref{eq:2.18}   show that 
$-\log (1- |r(x, k)|^2)/(1+k^2)$ is $L^1(\R)$. 
Since 
$r(x, k)/k$ is bounded, we deduce  from \eqref{eq:4.10} that $\re h(x, \cdot )$ 
is $L^1(\R)$.

When $\tau >0$, we  have
$$ h(x, \ii \tau) =\frac{1}{\pi} \int\limits_{\R} \tau \frac{h(x, k)}{k^2+\tau^2} \dd k =
 \frac{1}{\pi} \int\limits_{\R} \tau \frac{\re h(x, k)}{k^2+\tau^2} \dd k,
 $$
where we have used \eqref{eq:4.101} to get the second equality.
Hence
$$ \tau h(x, \ii \tau) = \frac{1}{\pi} \int\limits_{\R} \tau^2 \frac{\re h(x, k)}{k^2+\tau^2} \dd k
\to \frac{1}{\pi} \int \limits_{\R} \re h(x, k) \dd k \quad \text{when}\quad \tau\to \infty
 $$
We have thus obtained that 
$$ 
\begin{aligned}
 \int\limits_{\R} \frac{-\log (1- |r(x, k)|^2)}{k^2}\, \dd k &  =
 \int \limits_{\R} \re h(x, k) \dd k= \pi  \lim_{\tau\to \infty} \tau h(x, \ii \tau)\\
 & = \lim_{\tau\to \infty} \pi \int\limits_x^\infty q(y)r(y, \ii \tau) \dd y= \pi \int \limits_x^\infty Q^2(y) \dd y.
 \end{aligned}
 $$
 This proves \eqref{plancherel:1}.

We let $x\to -\infty$ in 
\eqref{eq:4.10} and use
\eqref{eq:3.9b}. 
It follows that
$$
\frac{-\log (1- |R_2(k)|^2)}{k^2} = 
-\frac{1}{2}\int\limits_{-\infty}^\infty \frac{q(y)(r(y, k)- \overline{r(y, k)})}{\ii k} \dd y .
$$
The equality \eqref{plancherel:2} follows then  by the same argument as before, with $h(x, k)$ replaced by $h(k) = -\int \frac{q(y) r(y, k)}{\ii k} \dd y$.

We see from \eqref{plancherel:2} that $ \int_{\R} Q^2(y) \dd y$ is uniquely determined by $R_2(k)$.  
Then the last statement in the theorem is a consequence  of the equality
$2 Q= Q^2+q$, and of the fact that  $\int\limits q(x)  \dd x$ is uniquely given 
by the limit of  $2\ii R(k)/k$ when $k\to 0$ 
(see lemma~\ref{lemma:2.2}).
\end{proof}

\section{The transmission coefficient}\label{sect:T}

The main result of this section shows  that  $R_2(k)$, 
$k\in \R$,
uniquely determines the transmission coefficient 
$T(k)$, $k\in \R$.

We have seen in Section~\ref{sect2}  that $T(\cdot)$ extends analytically to
$\Pi$,  continuously on  $\im k\ge 0$,  and that $\vert T(k)\vert \le 1$ when $k\in \R$.
We  first obtain estimates for  $T(k)$ when $k\in \overline{\Pi}$.
The main result in this direction is contained in the next proposition.

\begin{proposition}\label{prop:4.1}
Set $Q=1-c ^{-1}$ and denote
\begin{equation}\label{def:Ttilde}
 \tilde{T}(k)= T(k)\ee^{\ii k \int\limits_{\R} Q \dd x}.
\end{equation}
Then $
|\tilde{T}(k)|\le 1 $
when $k\in \overline{\Pi}$ .
\end{proposition}

To prove this proposition we need some lemmas.
We first establish some formulas for $\tilde{T}(k)$.

\begin{lemma}\label{lemma:Tform}
We have
\begin{equation} \label{eq:4.8}
%\begin{aligned} 
\tilde{T}(k)  
%&
= \lim_{x\to -\infty} \frac{2 
\ii k \ee^{\ii k \int\limits_{x}^\infty Q \dd x}}{m_1'(x, k) + 2\ii k m_1(x, k)} 
%\\ &
 = \exp\Big( \frac{\ii k}{2} \int\limits (Q^2 (y) - q(y) r(y, k)) \dd y\Big ) 
%\end{aligned}
\end{equation}
for every $k\in \overline{\Pi}$ .
\end{lemma}

\begin{proof}
From \eqref{eq:2.2} and \eqref{eq:m1} we have that 
$$
\frac{1}{T(k)} =\frac{1}{2\ii k}\Big [2\ii k -k^2\int\limits_{-\infty}^{\infty}q(t)m_1(t, k)\dd t \Big] 
=\lim_{x\to -\infty}\frac{(m_1'+2\ii km_1)(x, k)}{2\ii k}.
$$
Thus, if 
$ f(x, k):=2\ii k((m_1'+2\ii km_1)(x, k))^{-1}$, 
we have proved that
\begin{equation}\label{eq:4.6}
\tilde{T}(k)=
\lim_{x\to -\infty} f(x, k)\ee^{\ii k\int\limits_x^\infty  Q} = 
 \lim_{x\to -\infty} \frac{2 \ii k \ee^{\ii k \int\limits_{x}^\infty Q \dd x}}{m_1'(x, k) + 2\ii k m_1(x, k)}. 
\end{equation}
This proves the first equality in the statement.

Equation \eqref{eq:m1} yields
\begin{equation}\label{eq:4.7}
f' =-\frac{2\ii k (m_1''+ 2\ii k m_1')}{(m_1'+2\ii k m_1)^2}
=\frac{\ii k q}{2}(r+1) f,
\end{equation}
where we have used  the equality  $\displaystyle r+1=\frac{2\ii k m_1}{m'+2\ii k m_1 }$.
We also have $f(x, k)\to 1$ when $x\to \infty$, so we conclude  that 
$$f(x, k)=\ee^{-\ii k\int\limits_x^\infty \frac{q}{2}(r+1)}.
 $$
Thus, since  $q=2Q-Q^2$, we obtain the second equality in the statement.
\end{proof}

\begin{lemma}\label{lemma:4.2}
One has 
\begin{equation}\label{eq:4.2}
\big(2\kappa m_1 (x,\ii \kappa)-m_1 '(x, i\kappa)\big)
\ee^{\kappa \int\limits_x^\infty Q}\ge 2\kappa
\end{equation}
when $x\in \R$ and $\kappa\ge 0$.
\end{lemma}
   
\begin{proof}
The estimate \eqref{eq:4.2} is obvious for 
$\kappa=0$, both sides being equal to zero.

Fix $\kappa>0$. 
Then since $w(x, \ii \kappa)$ is real and  
$$
w(x, \ii \kappa)=\frac{(u_1'u_1)(x, \ii\kappa)}{-\kappa u_1^2(x, \ii \kappa)}>0,
$$
we get  $(u_1'u_1)(x, \ii \kappa)<0$. 
The function $x\to u_1 (x, \ii \kappa)$ does not vanish at any point of  $\R$ and is continuous,
therefore it has constant sign on  $\R$. 
Since $m_1 (x, \ii\kappa)=\ee^{\kappa x}u_1 (x, \ii \kappa)\to 1$ when $x\to \infty$, 
$u_1 (x, \ii \kappa)$ must be positive. 
Thus $u_1 '(x, \ii \kappa)<0$, and it follows that
$$
m_1'(x, \ii \kappa)<\kappa m_1(x, \ii \kappa). 
$$  

We have that
$$
u_1''(x, \ii \kappa)-\frac{\kappa^2}{c^2(x)}u_1 (x, \ii \kappa)=0.
$$
Therefore if we set $v(x)=u_1'(x, \ii \kappa) u_1 (x,\ii \kappa)$, $v$ satisfies
$$
v'=\frac{\kappa^2}{c^2}u_1 ^2(\cdot, \ii\kappa)+(u_1 ')^2(\cdot, \ii\kappa)\ge -2\frac{\kappa}{c}v.
$$
Hence  
$$ 
v(x) \ee^{2\kappa(x+\int\limits_x^\infty Q)}\le 
\lim_{x\to \infty} (u_1 'u_1 )(x, \ii \kappa)\ee^{2\kappa(x+\int\limits_x^\infty Q)}=-\kappa, 
$$ 
that is,
\begin{equation}\label{eq:4.4}
m_1(x, \ii \kappa)[\kappa m_1(x, \ii \kappa)-m_1'(x, \ii \kappa)] 
\ee^{2\kappa\int\limits_x^\infty Q}\ge \kappa.
\end{equation}
Since  $m_1(x, \ii\kappa)>0$, $\kappa m_1(x, \ii\kappa)-m_1'(x, \ii \kappa)>0$,
and  by using \eqref{eq:4.4}, we obtain
$$
\begin{aligned}
\frac{(2\kappa m_1-m_1')(x, \ii\kappa)}{2\kappa} & =  
\frac{(\kappa m_1-m_1')(x, \ii\kappa)}{2\kappa}+\frac{m_1(x, \ii \kappa)}{2}\\
&\ge \Big[ \frac{m_1(x, \ii \kappa)(\kappa m_1-m_1')(x, \ii \kappa)}
{\kappa}\Big]^{1/2}\ge  \ee^{-\kappa\int\limits_x^\infty Q}.
\end{aligned}
$$
This completes the proof.
\end{proof}

\begin{proof}[Proof of Proposition~\ref{prop:4.1}.]
From \eqref{eq:4.8} and the fact that $|r(x, k)|<1$ we see that there is an $a\ge 0$ such that
$\vert \tilde{T}(k)\vert \le \ee^{a|k|}$ when $\im k \ge 0$.
Thus $\tilde{T} $ is of angular order $1$ on $\{z \mid \im z >0, \re z>0\}$. 
We also have that $\vert \tilde{T}(k)\vert \le 1 $ when $k \in \R$ (from \eqref{eq:2.20}) and 
$k \in \ii \R$, 
from the first equality in \eqref{eq:4.8} and Lemma~\ref{lemma:4.2}.
The inequality  in the statement is then a consequence of
the  Phragm\'en-Lindel\"of
principle.
\end{proof}

The next theorem is the main result of this section.

\begin{theorem}\label{thm:4.3}
If  $c$  satisfies [H1] and [H2],
the transmission coefficient $(T(k))_{k\in\R}$ is uniquely determined 
by the reflection coefficient
$(R_2(k))_{k\in \R}$. 
\end{theorem}

\begin{proof} Recall that the function
$$
\tilde{T}(k)=T(k)\ee^{\ii k \int Q }, \qquad k\in \overline{\Pi}
$$
belongs to $H^\infty(\Pi)$ (Proposition~\ref{prop:4.1}), is continuous on  $\{\im k=0\}$ and has no zeros.
Note that $|\tilde{T}(k)|^2=|T(k)|^2=1- |R_2(k)|^2$ when $k$ is real, hence $|\tilde{T}(k)|$, $k\in \R$, 
is uniquely determined by $R_2(k)$.
Moreover, it follows from
\eqref{eq:4.8}  that
$$\tilde{T}(k) =
\exp\Big (-\ii k\int\limits_\R q(y)\frac{r(y, k)}{2}\ \dd y +
\ii k \int\limits_\R \frac{Q(y)^2}{2} \ \dd y\Big). $$
Lemma~\ref{lemma:4.5} shows that 
$$\lim_{\ka \to \infty} \frac{\log |\tilde{T}(\ii \ka)|}{\ka}=0.$$
It follows moreover that the mapping $k\to \log\vert\tilde{T}(k)\vert$ is harmonic, therefore 
$\tilde{T}$ is outer.
Then the factorization theorem of $H^\infty(\C)$ functions
gives that
\begin{equation}\label{eq:4.13}
\tilde{T}(z) =\gamma\Theta_F(z), \quad z\in \Pi, 
\end{equation}
where $\gamma\in \C$, $|\gamma|=1$, and 
$$
\Theta(z)=\exp\Big(-\frac{\ii }{2\pi }\int\limits_\R \big (\frac{1}{k-z}-\frac{k}{k^2+1}\Big)
\log |\tilde{T}(k)|^2\dd k \Big).
$$
Note that $\Theta$ is uniquely determined by the values  of $\tilde{T}$ on the real axis, hence by $R_2$.
It remains to determine $\gamma$.
 
From Lemma~\ref{lemma:Tform} we see that
$\log |\tilde{T}(k)|^2/k$ is locally $L^1$ and bounded. 
Then when $\tau\to 0$, 
$$\int\limits_\R \Big (\frac{1}{k-\ii \tau}-\frac{k}{k^2+1}\Big)
\log |\tilde{T}(k)|^2\dd k \to 
\int\limits_\R \frac{1}{k(k^2+1)}
\log |\tilde{T}(k)|^2\dd k=0$$
where in the last equality we have used that the integrand is an odd function. 
(Recall that  $\overline{T(k)}=T(-k)$ when $k$ is real.) 
Hence  $\Theta(\ii \tau )\to 1$ when $\tau\to 0$.
On the other hand $\tilde{T}(0)=T(0)=1$.
Hence 
$\gamma=1$, and this proves the theorem.
\end{proof}

\section{Further properties of $m_1(x, k)$}\label{sect:further}

Recall that when  $f$ is in the Nevanlinna class $N^+(\Pi)$ then the boundary values
$$ f(\lambda ) =\lim_{z\to \lambda} f(z)$$
exist non-tangentially almost everywhere $\lambda \in \R$, and $f$ can be recovered from these boundary values. 
Therefore one can identify $f$ with its boundary values at $\im z=0$.
From the maximum principle we also have that
$$L^2(\R)\cap N^+(\Pi)  =H^2(\Pi).$$  
(See \cite[Lemma~5.21]{RR}.)

Recall that we have set 
\begin{equation}\label{def:tildet}
\tilde{T}(k) = T(k) e^{\ii k\int Q}, \quad k \in \overline{\Pi}.
\end{equation}
Then $\tilde{T}$ is an $H^\infty(\Pi)$ outer function
(see the proof of theorem~\ref{thm:4.3}).
Denote
\begin{equation}\label{def:tildem}
\tilde{m}_1(x,k) = m_1(x, k) \ee^{-\ii k\int\limits_x^\infty Q}, \; \tilde{m}_2(x,k) = m_2(x, k) \ee^{-\ii k\int\limits^x_{-\infty} Q} \quad x\in \R, \, k \in \overline{\Pi}.
\end{equation}

\begin{lemma}\label{lemma:h2} Assume $x\in \R$ is fixed.
\begin{itemize}
 \item[a)]
 The function
$$\Pi\ni k\rightarrow \tilde{m}_1(x, k) \in \C
$$
belongs to $N^+(\Pi)$.
\item[b)] The function 
$$\Pi\ni k\rightarrow \frac{\tilde{T}(k)(\tilde{m}_1(x, k)-1)}{k}\in \C
$$
belongs to the Hardy space $H^2(\Pi)$.
\end{itemize} 
\end{lemma}

\begin{proof} Let $x$ be fixed. 
We first show that 
$k\to  \tilde{m}_1(x, k)$
belongs to the Nevanlinna class $N^{+}(\Pi)$.

For that, we note that if  $T_x(k)$ is the transmission coefficient associated to
$$c_x(y)=\begin{cases} c(y), & y\ge x\\
1, & y< x
\end{cases}
$$
then from \eqref{eq:2.12} we get
$$\frac{2\ii k}{T_x(k)}=u_1'(x, k)\ee^{-\ii k x} +\ii k \ee^{-\ii kx} u_1(x, k),
\qquad k\in \overline{\Pi},$$
hence
$$2 (T_x(k)u_1(x, k)e^{-\ii k x})^{-1} = w(x, k)+1.$$
It follows that
$$T_x(k)m_1(x, k)=r(x, k)+1 \quad \text{when $k\in \overline{\Pi}$},
$$
and thus
$$ \tilde{m}_1(x, k) = m_1(x, k)\ee^{-\ii k\int\limits_x^\infty Q}=
(1+r(x, k))(T_x(k)\ee^{\ii k\int\limits_x^\infty Q})^{-1}
$$
belongs to  $N^{+}(\Pi)$ since
$1+r(x, \cdot)$ is an $H^\infty(\Pi)$ function, while
$T_x(k)\ee^{-\ii k\int\limits_x^\infty Q}$ is an 
$H^\infty(\Pi)$ outer function.

Since $T(k)\ee^{-\ii k\int Q}$ is $H^\infty (\Pi)$ and outer we get that 
$f(k)=\tilde{T}(k)
m_1(x, k)\ee^{-\ii \ka \int\limits_x^\infty Q}$ belongs  ${N}^{+}(\Pi)$.

By using the maximum principle and   
\eqref{eq:2.4i}
it suffices to show that $f/(\ii +k)$ is $L^2$ on $\R$.
We note that  lemma~\ref{lemma:2.3} 
leads to
$$ |f(k)/(\ii +k)|^2 \le 2 \re (T(k) m_1 (x, k) m_2(x, k))/(1+k^2).$$
On the other hand it follows from \eqref{eq:2.12}  that
$$
\begin{aligned}
\re (T(k) m_1(x, k) m_2(x, k))& =  \re (\tilde{T}(k) \tilde{m}_1(x, k)\tilde{m}_2(x, k))\\
& =\re \Big(\frac{1}{w(x, k)+ w_-(x, k)}\Big),
\end{aligned}
$$
and this, along with Corollary~\ref{cor:poz},  shows that the function   
$$k \to \re (T(k) m_1 (x, k) m_2(x, k))/(1+k^2)$$
 belongs to  $L^1$.
We have obtained  that $f/(\ii +k)$ belongs to  $L^2(\R)$, which completes the 
proof.
\end{proof}

\begin{lemma}\label{growth}
 Assume $x\in \R$ is fixed. 
Then  
the functions
$$ \R\ni k \mapsto \frac{\vert \tilde{T}(k)\tilde{m}_j(x, k)\vert^2}{1+k^2}\in \R, \quad j=1, 2, $$ 
are in $L^1(\R)$. 
\end{lemma}

\begin{proof}
We prove the lemma for $\tilde{m}_1$. 
The proof is completely similar for  $\tilde{m}_2$. 
We use the notation in Lemma~\ref{lemma:h2}. 
By Lemma~\ref{lemma:2.3} and the definition of $\tilde{T}$ and $\tilde{m}_1$ we have that 
$$ \vert \tilde{T}(k)\tilde{m}_1(x, k)\vert^2 \le 2  \re (\tilde{T}(k) \tilde{m}_1(x, k)\tilde{m}_2(x, k))=
\re \Big(\frac{1}{w(x, k)+ w_-(x, k)}\Big).$$
Now the lemma follows from the fact that $w$ and $w_-$ have positive real parts on $\overline{\Pi}$ and 
are bounded on $\ii \R$. 
\end{proof}

\begin{lemma}\label{atinfty}
For every $x$ and every $\epsilon>0$ there is a constant $C_\epsilon>0$ such that
$$ |\tilde{m}_j(x, \ii \kappa))|\le C_\epsilon \ee^{\epsilon\kappa}\; \text{for all} \;  \kappa>0,  $$
 $j=1, 2$. The constant $C_\epsilon$ may be chosen independently on $x$. 
\end{lemma}

\begin{proof}
We prove the lemma for $\tilde{m}_1$. 
The proof  for $\tilde{m}_2$ is completely similar.
Note first that 
$\tilde{m}_j(x, \ii \kappa)$ is positive when $\kappa>0$. 
We use the notation in the proof of Lemma~\ref{lemma:h2}. Recall that
$$ \tilde{m}_1(x, k)= (1+r(x, k))(T_x(k)\ee^{\ii\int\limits_{x}^\infty Q})^{-1}.$$ 
Hence 
$$ 
\begin{aligned}
\frac{\log(\tilde{m}_j(x, \ii \kappa))}{\kappa} & = 
\frac{\log(1+r(x, \ii\kappa))}{\kappa} -
 \frac{\log(T_x(k)\ee^{\ii\int\limits_{x}^\infty Q})}{\kappa}   \\
& = 
\frac{\log(1+r(x, \ii\kappa))}{\kappa} +  \frac{1}{2} \int\limits_x^\infty ( Q^2(y) - q(y) r(y, k))\dd y.
\end{aligned} 
$$
We have used here \eqref{eq:4.8} written for $T_x$. 
Then the lemma follows from the fact that $\vert r(x, k)\vert \le 1$ on $\overline{\Pi}$ and by Lemma~\ref{lemma:4.5}.
\end{proof}

\section{The uniqueness}\label{sect:uniq}

In this section we prove the main result, Theorem~\ref{thm:4.1.1}.

\begin{proof}[Proof of theorem~\ref{thm:4.1.1}]
First set $R_2(k):=R_2(k, c_1)= R_2(k, c_2)$, and   note that by Theorem~\ref{thm:4.3} and 
Lemma \ref{lemma:4.6} one also has that 
$$
T(k; c_1)=T(k, c_2)=:T(k)\ \text{and}\ 
\int\limits_\R Q_1 \dd s=\int\limits_\R Q_2\dd s.$$

We denote by $u_{1, 1}(x, k)$ and $u_{1,2}(x, k)$ the Jost solutions corresponding to
 $c_1$, 
and by $u_{2, 1}(x, k)$ and $u_{2,2} (x, k)$ those corresponding to $c_2$. 
We write
$$ u_{1, j}(x, k) =\ee^{\ii kx} m_{1, j}(x, k), \ u_{2, j}(x, k)=\ee^{-\ii kx} m_{2, j}(x, k),\,  j=1, 2.$$
With this notation \eqref{eq:2.11ii} becomes
\begin{equation}\label{eq:4001}
T(k) m_{1, j}(x, k) =R_2(k) \ee^{-2\ii kx} m_{2, j}(x, k) +\overline{m_{2, j}(x, k)}, \ j=1, 2
\end{equation}
when $k$ is real.
For $k\in \overline{\Pi}$ we define
$$
\begin{aligned}
& \tilde{T}(k) =\ee^{\ii k\int Q_1} T(k) =
\ee^{\ii k\int Q_2} T(k), \\ 
& \tilde{m}_{1, j}(x, k)=\ee^{-\ii k\int\limits_x^\infty  Q_j} m_{1, j}(x, k), \,  
\tilde{m}_{2, j}(x, k)=\ee^{-\ii k\int\limits_{-\infty}^x  Q_j} m_{2, j}(x, k).  
\end{aligned}
$$ 
Note that from Lemma~\ref{atinfty} we have that for fixed $x$ and  for every $\ep>0$ there is a constant $C_\ep$, independnt of $x$,  
such that 
\begin{equation}\label{compact:1}
 \abs{\tilde{m}_{1, j}(x, \ii \kappa)}\le C_\ep \ee^{\ep \kappa},\; \abs{\tilde{m}_{2, j}(x, \ii \kappa)}\le C_\ep \ee^{\ep \kappa},  
\end{equation}
when $\kappa> 0$. 

Then by  \eqref{eq:4001} we  get 
\begin{equation}\label{eq:4002}
\tilde{T}(k) \tilde{m}_{1, j}(x, k)=
R_2(k) \ee^{-2\ii k (x-\int\limits_{-\infty}^x Q_j)} \tilde{m}_{2, j}(x, k)
+\overline{\tilde{m}_{2, j}(x, k)}, \ j=1, 2.
\end{equation}
We now change variables 
$$ y:= \chi_j(x) = x-\int\limits_{-\infty}^x Q_j, $$
and note that $\chi_j'(x) = 1/c_j(x)>0$. 
In the new variables, \eqref{eq:4002} becomes 
\begin{equation}\label{compact:2}
\tilde{T}(k)  \tilde{m}_{1, j}(\chi_j^{-1}(y), k)=
R_2(k) \ee^{-2\ii k y} \tilde{m}_{2, j}(\chi_j^{-1}(y), k)
+\overline{\tilde{m}_{2, j}(\chi_j^{-1}(y), k)},
\end{equation}
whenever $k \in \R$ and for $j=1, 2$.  
Set now
\begin{equation}\label{def:v}
\begin{aligned}
v_{1}(y, k) & =  \ee^{\ii k y} (\tilde{m}_{1, 1}(\chi_1^{-1}(y), k) -   \tilde{m}_{1, 2}(\chi_2^{-1}(y), k)),      \\
v_{2}(y, k) &  = \ee^{-\ii k y}(\tilde{m}_{2, 1}(\chi_1^{-1}(y), k)- \tilde{m}_{2, 2}(\chi_2^{-1}(y), k)),  
\end{aligned}
\end{equation}
for all $k\in \overline{\Pi}$ and $y\in \R$.  
Thus \eqref{compact:2} becomes 
\begin{equation}\label{eq:4002.1.1}
\tilde{T}(k) v_{1}(y, k)  =
R_2(k) v_{2}(y, k)
+\overline{v_{2}(y, k)},  
\end{equation}
for all $k\in \R$ and $y\in \R$.

Let now $x >y$ be fixed. 
Denote $g(x, y,  k)= \tilde{T}(k) v_1(x, k) v_2(y, k)$.
 Then by Lemma~\ref{slemma} and \eqref{eq:4002.1.1}  it follows that 
 \begin{equation}\label{3001}
 \begin{aligned}
 2\re g(x, y, k) & = \abs{\tilde{T}(k)}^2 v_1(x, k) \overline{v_1(y, k)}+ \abs{\tilde{T}(k)}^2 v_2(x, k) 
 \overline{v_2(y, k)},\\
  & = \re(\abs{\tilde{T}(k)}^2  v_1(x, k) \overline{v_1(y, k)}+ \abs{\tilde{T}(k)}^2 v_2(x, k) 
 \overline{v_2(y, k)}),   
  \quad k\in \R.
  \end{aligned}
\end{equation}
Then  $k\to \re g(x, y, k)/k^2$ is a continuous, even,  $L^1(\R)$ function in variable $k$. 
 Also, $\Pi\ni k \to g(x, y, k)$ belongs to the Nevanlinna class, 
 since 
 \begin{equation}\label{g:Nev}
 g(x, y, k) = \ee^{\ii k(x- y)}  h(x, y, k), 
 \end{equation}
 where
 $$ \begin{aligned}
h(x, y, k) =  \tilde{T}(k) & (\tilde{m}_{1, 1}(\chi_j^{-1}(x), k) -   \tilde{m}_{1, 2}(\chi_2^{-1}(x), k)) \times \\\times 
& (\tilde{m}_{2, 1}(\chi_1^{-1}(y), k) -   \tilde{m}_{2, 2}(\chi_2^{-1}(y), k))
\end{aligned} 
$$
 and all factors in the right hand-side of \eqref{g:Nev} are Nevanlinna (Lemma~\ref{lemma:h2},
 Prop.~\ref{prop:4.1}).  
The function  $g(x, y, \cdot )$
satisfy the conditions of Corollary~\ref{cor1001}. 
Indeed $g(x, y, \cdot)$ is bounded on $\ii \R$, by Lemma~\ref{atinfty}.
Moreover $\overline{g(x, y, k)}= g(x, y, -\overline{k})$ when $k\in \overline{\Pi}$, $g(x, y, 0)=0$
and  $k\to \re{g(x, y, k)}/k^2$ is $L^1(\R)$.  
Also \eqref{3001} and  Lemma~\ref{growth} 
ensures that the condition \eqref{cor1001_3}  in Corollary~\ref{cor1001} is satisfied. 
Then it follows that 
 \begin{equation}\label{imp1}
 g(x, y, k)= 
 \frac{1}{\pi \ii} \int \re g(x, y, t) \Big (\frac{1}{t-k} -\frac{t}{t^2+1}\Big) \, \dd t , \quad k \in \Pi. 
 \end{equation} 
 
We denote by $g(x, k) = g(x, y, k)$. 
Then in \eqref{imp1} we can let $x\to y$ and obtain that 
 \begin{equation}\label{imp2}
 g(x, k)= 
 \frac{1}{\pi \ii} \int \re g(x, t) \Big (\frac{1}{t-k} -\frac{t}{t^2+1}\Big) \, \dd t , \quad k \in \Pi. 
 \end{equation} 
Note that 
\begin{equation}\label{imp25}
 \re g(x, k)= \abs{\tilde{T}(k)}^2( \abs{v_1(x, k)}^2+ \abs{v_2(x, k)}^2) \ge 0, \quad k \in \R. 
\end{equation}
Take now $k=\ii \kappa$ in \eqref{imp2}, and obtain
\begin{equation}\label{imp3}
 \frac{g(x, \ii \kappa)}{ \kappa} = 
 \frac{1}{\pi} \int \re g(x, t )\frac{1}{t^2+\kappa^2}  \, \dd t , \quad \kappa>0. 
\end{equation}
Now if we let $\kappa\to 0$ in the left-hand side of \eqref{imp3}, since  $g(x, \ii \kappa)/ \kappa$ converges to $0$, 
and $t\to g(x, t)/t^2$ is $L^1$, we get that 
$$  \frac{1}{\pi} \int \frac{\re g(x, t )}{t^2}  \, \dd t= 0.$$
This equality, the above \eqref{imp25} and the definition \eqref{def:v} show that  
$$  
\begin{aligned}
\tilde{m}_{1, 1}(\chi_1^{-1}(y), k) & = & \tilde{m}_{1, 2}(\chi_2^{-1}(y), k),      \\
\tilde{m}_{2, 1}(\chi_1^{-1}(y), k) & = & \tilde{m}_{2, 2}(\chi_2^{-1}(y), k) 
\end{aligned}
$$
for all $k\in \R$. 
Then from the limits of $\tilde{m}_{1, j}(\chi_j^{-1}(y), k)/k$ when $k\to 0$ we get that 
$$ \int\limits^{\chi_1^{-1}(y)} _{-\infty} Q_1(s) \, \dd s = \int\limits^{\chi_2^{-1}(y)}_{-\infty} Q_2(s) \, \dd s, $$
that is, 
$$ \chi_1^{-1}(y) - \chi_1( \chi_1^{-1}(y))  = 
\chi_2^{-1}(y) - \chi_2( \chi_2^{-1}(y)).$$ 
This shows that $\chi_1^{-1}(y) = \chi_2^{-1}(y)$ 
for every $y$, which in turn implies that $c_1= c_2$ a.e.
This finishes the proof. 
\end{proof}

\subsection*{Acknowledgments.} 
We thank the referee for the valuable comments and suggestions that considerably improved our paper. 

This research has been partially supported by the 
Laboratoire Europ\'een Associ\'e CNRS Franco-Roumain "Math-Mode".
The first author is also partially supported 
 by the Grant
of the Romanian National Authority for Scientific Research, CNCS-UEFISCDI,
project number PN-II-ID-PCE-2011-3-0131.

\enddocument